\documentclass[amsthm]{elsart}
\usepackage{euscript,amsmath, amssymb, amsfonts}

\newcommand{\R}{{\mathbb R}}
\newcommand{\C}{{\mathbb C}}
\newcommand{\di}{{\rm d}}
\newcommand{\supp}{\mathrm{supp}\,}

\begin{document}

\begin{frontmatter}

\title{On localization properties of Fourier transforms of hyperfunctions}
\author{A.~G.~Smirnov\thanksref{th}}
\ead{smirnov@lpi.ru}
\address{I.~E.~Tamm Theory Department, P.~N.~Lebedev
Physical Institute, Leninsky prospect 53, Moscow 119991, Russia}

\thanks[th]{The research was supported by the grant LSS-1615.2008.2.}

\begin{abstract}
In [Adv.~Math. {\bf 196} (2005) 310--345] the author introduced a new generalized function space $\mathcal U(\mathbb R^k)$ which can be naturally
interpreted as the Fourier transform of the space of Sato's hyperfunctions on $\mathbb R^k$. It was shown that all Gelfand--Shilov spaces $S^{\prime
0}_\alpha(\mathbb R^k)$ ($\alpha>1$) of analytic functionals are canonically embedded in $\mathcal U(\mathbb R^k)$. While the usual definition of
support of a generalized function is inapplicable to elements of $S^{\prime 0}_\alpha(\mathbb R^k)$ and $\mathcal U(\mathbb R^k)$, their localization
properties can be consistently described using the concept of {\it carrier cone} introduced by Soloviev [Lett.~Math.~Phys. {\bf 33} (1995) 49--59;
Comm.~Math.~Phys. {\bf 184} (1997) 579--596]. In this paper, the relation between carrier cones of elements of $S^{\prime 0}_\alpha(\mathbb R^k)$ and
$\mathcal U(\mathbb R^k)$ is studied. It is proved that an analytic functional $u\in S^{\prime 0}_\alpha(\mathbb R^k)$ is carried by a cone
$K\subset\mathbb R^k$ if and only if its canonical image in $\mathcal U(\mathbb R^k)$ is carried by $K$.
\end{abstract}

\end{frontmatter}

\section{Introduction}
\label{s1}

It is well known from the theory of (Fourier) hyperfunctions that the description of the localization properties of generalized functions becomes a
nontrivial problem in the case when all test functions are analytic (see, e.g., Chap.~9 in~\cite{Hoermander}). The standard definition of support is
inapplicable to such generalized functions because of the lack of test functions with compact support. In particular, this difficulty arises in the
case of the Gelfand--Shilov spaces $S^\beta_\alpha(\R^k)$ with $0\leq \beta<1$, which consist of (restrictions to $\R^k$ of) entire analytic
functions on $\C^k$ (we refer the reader to \cite{GS} for the definition and basic properties of Gelfand--Shilov spaces). It was shown by
Soloviev~\cite{Sol1,Sol2} that the localization properties of elements of $S^{\prime \beta}_\alpha(\mathbb R^k)$ (topological dual of
$S^\beta_\alpha(\mathbb R^k)$) can be consistently described using the concept of {\it carrier cone} instead of support. The definition of carrier
cones is based on introducing, for every closed cone $K$, a suitable test function space $S^\beta_\alpha(K)$ in which $S^{\beta}_\alpha(\R^k)$ is
densely embedded (the precise definition will be given later in this section); a functional $u\in S^{\prime \beta}_\alpha(\R^k)$ is said to be
carried by a closed cone $K$ if $u$ has a continuous extension to $S^\beta_\alpha(K)$. Functionals carried by a closed cone $K$ have much the same
properties as the ordinary generalized functions whose support is contained in $K$. In particular, every element of $S^{\prime \beta}_\alpha(\R^k)$
has a unique minimal carrier cone~\cite{Sol1}.

In~\cite{AIM}, we introduced a new generalized function space $\mathcal U(\mathbb R^k)$ which can be naturally interpreted as the Fourier transform
of the space of Sato's hyperfunctions on $\mathbb R^k$. The space of hyperfunctions on $\R^k$ can be thought of as the limiting case as
$\alpha\downarrow 1$ of the ultradistribution spaces $S^{\prime \alpha}_0(\R^k)$. Therefore, it is natural to try to define the Fourier transform
$\mathcal U(\mathbb R^k)$ of the space of hyperfunctions by passing to the limit $\alpha\downarrow 1$ in the definition of the spaces $S^{\prime
0}_\alpha(\R^k)$, which are the Fourier transforms of $S^{\prime \alpha}_0(\R^k)$ (recall that the Fourier transformation just interchanges the
indices of Gelfand--Shilov spaces). However, we cannot just set $\mathcal U(\mathbb R^k)=S^{\prime 0}_1(\R^k)$ because the space $S^0_1(\R^k)$ is
trivial. In~\cite{AIM}, we proposed a procedure for making $S^{\prime 0}_1(\R^k)$ into a nontrivial space. The key observation is that the spaces
$S^0_\alpha(K)$ over {\it proper}\footnote{A cone $U$ in $\R^k$ will be called proper if $\bar U\setminus\{0\}$ is contained in an open half-space of
$\R^k$ (the bar denotes closure). For convex closed cones, this definition is equivalent to the usual one, by which a cone is called proper if it
contains no straight lines. } cones remain nontrivial after passing to the limit $\alpha\downarrow 1$. This allows us to construct $\mathcal
U(\mathbb R^k)$ by suitably ``gluing together'' the generalized function spaces $S^{\prime 0}_1(K)$ associated with proper closed cones $K\subset
\R^k$ (the precise meaning of such gluing is given by Definition~\ref{dxx3}).

The properties of the elements of $\mathcal U(\R^k)$, which we called \emph{ultrafunctionals}, are quite similar to those of analytic functionals in
$S^{\prime 0}_\alpha(\R^k)$. In particular, the definition of carrier cones is extended to the case of the space $\mathcal U(\R^k)$ and it turns out
that every ultrafunctional has a uniquely defined minimal carrier cone. Moreover, for any $\alpha>1$, there is a natural mapping $S^{\prime
0}_\alpha(\R^k)\to \mathcal U(\R^k)$. The aim of this paper is to prove the following relation between carrier cones of elements of $S^{\prime
0}_\alpha(\mathbb R^k)$ and~$\mathcal U(\mathbb R^k)$.

\begin{thm}
\label{t1} Let $\alpha>1$. The canonical mapping $\varepsilon^\alpha\colon S^{\prime 0}_\alpha(\R^k) \to \mathcal U(\R^k)$ is injective. A nonempty
closed cone $K$ is a carrier cone of a functional $u\in S^{\prime 0}_\alpha(\R^k)$ if and only if $K$ is a carrier cone of $\varepsilon^\alpha u$.
\end{thm}

The injectivity of $\varepsilon^\alpha$ means that $S^{\prime 0}_\alpha(\mathbb R^k)$ can be considered as a subspace of $\mathcal U(\mathbb R^k)$.
In fact, it has already been established in~\cite{AIM} using the injectivity of the canonical mappings of the spaces of ultradistributions to the
space of hyperfunctions. In this paper, however, we shall give a direct proof of the injectivity of $\varepsilon^\alpha$ that does not appeal to the
properties of hyperfunctions.

Before we pass to the proof of this theorem, we first need to give precise definitions of carrier cones, ultrafunctionals, and the mapping
$\varepsilon^\alpha$. As mentioned above, carrier cones can be consistently defined for all spaces $S^\beta_\alpha$ with $0\leq \beta<1$, but we
shall confine ourselves to the spaces $S^0_\alpha$ entering the formulation of Theorem~\ref{t1}. Throughout the paper, all cones are supposed to be
nonempty. We say that a cone $W$ is a conic neighborhood of a cone $U$ if $W$ contains $U$ and $W\setminus\{0\}$ is an open set (note that the
degenerate cone $\{0\}$ is a conic neighborhood of itself).

\begin{defn}
\label{d1} Let $\alpha\geq 1$ and $U$ be a cone in $\R^k$. The Banach space $S^{0,B}_{\alpha,A}(U)$ consists of entire analytic functions on $\C^k$
with the finite norm
\[
\|f\|^\alpha_{U,A,B}=\sup_{z\in \C^k} |f(z)|e^{-\sigma^\alpha_{U,A,B}(z)},
\]
where
\begin{equation}\label{000}
\sigma^\alpha_{U,A,B}(x+iy)= -|x/A|^{1/\alpha}+\delta_U(Bx)+|By|,
\end{equation}
$\delta_U(x)=\inf_{x'\in U}|x-x'|$ is the distance from $x$ to $U$, and $|\cdot|$ is a norm on $\R^k$. The space $S^0_\alpha(U)$ is defined by the
relation $S^0_\alpha(U)=\bigcup_{A,B>0,\,W\supset U} S^{0,B}_{\alpha,A}(W)$, where $W$ runs over all conic neighborhoods of $U$ and the union is
endowed with the inductive limit topology.
\end{defn}

Clearly, the definition of $S^0_\alpha(U)$ does not depend on the choice of the norm on $\R^k$. For definiteness, we assume the norm $|\cdot|$ to be
uniform, i.e.,  $|x|=\sup_{1\leq j\leq k}|x_j|$. If $U=\R^k$, then this definition is equivalent to the original definition of $S^0_\alpha(\R^k)$ due
to Gelfand and Shilov (see~\cite{GS}, Sec.~IV.2.3). If $U\subset U'$, then we obviously have the continuous inclusion $S^0_\alpha(U')\to
S^0_\alpha(U)$. Let $\rho^\alpha_{U,U'}$ denote the natural mapping from $S^{\prime 0}_\alpha(U)$ to $S^{\prime 0}_\alpha(U')$ (if $u\in S^{\prime
0}_\alpha(U)$ then $\rho^\alpha_{U,U'}u$ is the restriction of $u$ to $S^{0}_\alpha(U')$). For any $\alpha\geq 1$ and any cone $U\subset \R^k$,
$S^0_\alpha(U)$ is a nuclear DFS\footnote{Recall~\cite{Komatsu} that DFS spaces are, by definition, inductive limits of countable sequences of
locally convex spaces with compact linking mappings.} space (see~\cite{Smirnov}, Lemma~4).

Let $\alpha>1$. A closed cone $K$ is said to be a carrier cone of a functional $u\in S^{\prime 0}_\alpha(\R^k)$ if $u$ has a continuous extension to
the space $S^{0}_\alpha(K)$. The following basic properties of carrier cones were established in~\cite{Sol1,Sol2}.

\begin{thm} \label{t2}
Let $\alpha>1$. Then we have
\begin{enumerate}
\item[$1$.] The space $S^0_{\alpha}(\R^k)$ is dense in $S^0_{\alpha}(U)$ for any cone $U\subset\R^k$.

\item[$2$.] If both $K_1$ and $K_2$ are carrier cones of $u\in S^{\prime 0}_{\alpha}(\R^k)$, then so is $K_1\cap K_2$.

\item[$3$.] If $K_1$ and $K_2$ are closed cones in $\R^k$, then for any $u\in S^{\prime 0}_{\alpha}(\R^k)$ carried by $K_1\cup K_2$, there exist
$u_{1,2}\in S^{\prime 0}_{\alpha}(\R^k)$ carried by $K_{1,2}$ such that $u=u_1+u_2$.

\end{enumerate}
\end{thm}

Statement~1 in Theorem~\ref{t2} shows that the space of the functionals with the carrier cone $K$ is naturally identified with $S^{\prime
0}_{\alpha}(K)$ and that all mappings $\rho^\alpha_{U,U'}$ are injective, while statement~2 in Theorem~\ref{t2} implies that every functional in
$S^{\prime 0}_{\alpha}(\R^k)$ has a uniquely defined minimal carrier cone. The next result, which follows from Lemma~2.8 in~\cite{AIM}, shows that
$S^{\prime 0}_\alpha(\R^k)$ can be expressed in terms of the spaces $S^{\prime 0}_\alpha(K)$ over proper cones $K$.

\begin{lem}\label{l1}
Let $\alpha>1$, $K$ be a closed cone in $\R^k$, and $\mathcal P(K)$ be the set of all nonempty closed proper subcones of $K$ ordered by inclusion.
Then the embeddings $\rho^\alpha_{M,K}\colon S^{\prime 0}_\alpha(M)\to S^{\prime 0}_\alpha(K)$, $M\in\mathcal P(K)$, induce a topological isomorphism
\begin{equation}\label{1}
\varinjlim\nolimits_{M\in\mathcal P(K)} S^{\prime
0}_\alpha(M) \simeq S^{\prime 0}_\alpha(K),
\end{equation}
where the inductive limit is taken with respect to the linking mappings $\rho^\alpha_{M,M'}$, where $M,M'\in \mathcal P(K)$ are such that $M\subset
M'$.
\end{lem}

If $K$ is a proper closed cone, then the space $S^{0}_1(K)$ is nontrivial because it contains all exponentials $x+iy\to e^{l(x)+il(y)}$, where $l$ is
a linear functional on $\R^k$ such that $l(x)<0$ for any $x\in K\setminus\{0\}$. Lemma~\ref{l1} suggests that we can try to define the
``nontrivialization'' $\mathcal U(\R^k)$ of the space $S^{\prime 0}_1(\R^k)$ (and, more generally, of the space $S^{\prime 0}_1(K)$ over an arbitrary
closed cone $K$) as the left-hand side of~(\ref{1}) with $\alpha=1$. We then arrive at the following definition.

\begin{defn}
\label{dxx3} Let $K$ be a closed cone in $\R^k$. The space $\mathcal U(K)$ is defined as the inductive limit $\varinjlim\nolimits_{M\in \mathcal
P(K)} S^{\prime 0}_1(M)$, where $\mathcal P(K)$ is the set of all nonempty proper closed cones contained in $K$ and the inductive limit is taken with
respect to the linking mappings $\rho^1_{M,\,M'}$, $M,M'\in \mathcal P(K)$. The elements of $\mathcal U(\R^k)$ are called ultrafunctionals. A closed
cone $K$ is said to be a carrier cone of an ultrafunctional $u$ if the latter belongs to the image of the canonical mapping from $\mathcal U(K)$ to
$\mathcal U(\R^k)$.
\end{defn}

The canonical mapping $\mathcal U(K)\to \mathcal U(\R^k)$ used in Definition~\ref{dxx3} is induced by the mappings $\rho^1_{M,\,M'}$, where
$M\in\mathcal P(K)$ and $M'\in\mathcal P(\R^k)$ are such that $M\subset M'$. More generally, for any closed cones $K\subset K'$, the natural mappings
$\rho^1_{M,\,M'}$, where $M\in\mathcal P(K)$ and $M'\in\mathcal P(K')$ are such that $M\subset M'$, induce a canonical mapping $\rho^{\mathcal
U}_{K,\,K'}\colon {\mathcal U}(K)\to{\mathcal U}(K')$. In~\cite{AIM}, the following analogue of Theorem~\ref{t2} for ultrafunctionals was proved.

\begin{thm} \label{t3}
$\mbox{}$
\begin{enumerate}
\item[$1$.] The natural mapping $\rho^{\mathcal U}_{K,\,K'}\colon {\mathcal U}(K)\to{\mathcal U}(K')$ is injective for any closed cones $K$ and $K'$
such that $K\subset K'$.

\item[$2$.] Let $\{K_\omega\}_{\omega\in\Omega}$ be an arbitrary family of carrier cones of an ultrafunctional $u$. Then $\bigcap_{\omega\in\Omega}
K_{\omega}$ is also a carrier cone of $u$.

\item[$3$.] Let $K_1$ and $K_2$ be closed cones in $\R^k$ and an ultrafunctional $u$ be carried by $K_1\cup K_2$. Then there are $u_{1,2}\in
{\mathcal U}(\R^k)$ carried by $K_{1,2}$ such that $u=u_1+u_2$.

\end{enumerate}
\end{thm}

It follows from Statement~2 of Theorem~\ref{t3} that every ultrafunctional has a uniquely determined minimal carrier cone. For any proper closed cone
$K$, the space $\mathcal U(K)$ is naturally isomorphic to $S^{\prime 0}_1(K)$ and, therefore is nontrivial. In view of Statement~1 of
Theorem~\ref{t3}, this implies the nontriviality of $\mathcal U(K)$ for any closed cone $K$. If $\alpha>1$, then we have continuous inclusions
$S^0_1(K)\subset S^0_\alpha(K)$ for all closed cones $K$. The natural mappings $j^\alpha_K\colon S^{\prime 0}_\alpha(K)\to S^{\prime 0}_1(K)$ taking
an element of $S^{\prime 0}_\alpha(K)$ to its restriction to $S^0_1(K)$ are compatible with the linking mappings $\rho^1_{K,K'}$ and
$\rho^\alpha_{K,K'}$, i.e., $j^\alpha_{K'} \rho^\alpha_{K,K'}=\rho^1_{K,K'} j^\alpha_{K}$ for any $K\subset K'$. They therefore uniquely determine a
natural mapping from $\varinjlim_{M\in\mathcal P(K)}S^{\prime 0}_\alpha(M)$ to $\varinjlim_{M\in\mathcal P(K)}S^{\prime 0}_1(M)=\mathcal U(K)$ for
any closed cone $K$. Let $e^\alpha_K\colon S^{\prime 0}_\alpha(K)\to\mathcal U(K)$ be the composition of this mapping with the canonical isomorphism
$S^{\prime 0}_\alpha(K)\to \varinjlim_{M\in\mathcal P(K)}S^{\prime 0}_\alpha(M)$, which exists by Lemma~\ref{l1}. The mapping
$\varepsilon^\alpha\colon S^{\prime 0}_\alpha(\R^k)\to \mathcal U(\R^k)$ entering the formulation of Theorem~\ref{t1} is defined by setting
$\varepsilon^\alpha=e^\alpha_{\R^k}$.

The mappings $e^\alpha_K$ are compatible with the linking mappings $\rho^{\mathcal U}_{K,K'}$ and $\rho^\alpha_{K,K'}$:
\begin{equation}\label{compat}
e^\alpha_{K'}\, \rho^\alpha_{K,K'}=\rho^{\mathcal U}_{K,K'}\, e^\alpha_{K},\quad K\subset K'.
\end{equation}
If $u\in S^{\prime 0}_\alpha(\R^k)$ is carried by $K$, then we have $u=\rho^\alpha_{K,\R^k} v$ for some $v\in S^{\prime 0}_\alpha(K)$.
By~(\ref{compat}), this implies that $\varepsilon^\alpha u = \rho^{\mathcal U}_{K,\R^k} e^\alpha_K v$, i.e., $\varepsilon^\alpha u$ is carried by
$K$. Thus, it is only the ``if'' part of Theorem~\ref{t1} that needs proof.

Theorem~\ref{t1} is a generalization of Theorem~4 in~\cite{Sol2004} which states that, given $1<\alpha'<\alpha$, a closed cone $K$ is a carrier cone
of $u\in S^{\prime 0}_\alpha(\R^k)$ if and only if it is a carrier cone of the restriction of $u$ to $S^{0}_{\alpha'}(\R^k)$. Indeed, let $v$ be the
restriction of $u$ to $S^{0}_{\alpha'}(\R^k)$. Then we have $\varepsilon^\alpha u = \varepsilon^{\alpha'} v$ and the above statement follows
immediately from Theorem~\ref{t1}. It should be noted, however, that despite the similarity of formulations, the proof of Theorem~\ref{t1} turns out
to be considerably more complicated than that of Theorem~4 in~\cite{Sol2004}.

We shall prove Theorem~\ref{t1} in two stages. In the next section, we prove the statement under the additional assumption that $u$ is carried by
some {\it proper} closed cone containing $K$ (see Theorem~\ref{t4} below). While the treatment in Section~\ref{s2} is mostly analytic, passing to the
general case turns out to be a purely algebraic problem, which is solved in Section~\ref{s3} using the corresponding technique developed
in~\cite{AIM}.

\section{The case of proper cones}
\label{s2}

The aim of this section is to prove the next statement.

\begin{thm}\label{t4}
Let $\alpha\geq 1$ and $K\subset K'$ be proper closed cones in $\R^k$. If $u_1\in S^{\prime 0}_1(K)$ and $u_2\in S^{\prime 0}_\alpha(K')$ coincide on
$S^0_1(K')$, then there is $u\in S^{\prime 0}_\alpha(K)$ such that $u_1$ and $u_2$ are the restrictions of $u$ to $S^0_1(K)$ and $S^0_\alpha(K')$
respectively.
\end{thm}

Theorem~\ref{t4} implies Theorem~\ref{t1} under the assumption that $u$ is carried by some proper cone $K'\supset K$. Indeed, suppose
$\varepsilon^\alpha u$ is carried by $K$. Let $u_1\in\mathcal U(K)$ and $u_2\in S^{\prime 0}_\alpha(K')$ be such that $\varepsilon^\alpha
u=\rho^{\mathcal U}_{K,\R^k}u_1$ and $u=\rho^\alpha_{K',\R^k}u_2$. Using~(\ref{compat}), we obtain $\rho^{\mathcal U}_{K',\R^k} e^\alpha_{K'} u_2 =
\rho^{\mathcal U}_{K,\R^k}u_1=\rho^{\mathcal U}_{K',\R^k}\rho^{\mathcal U}_{K,K'}u_1$. By Statement~1 of Theorem~\ref{t3}, this implies that
$\rho^{\mathcal U}_{K,K'}u_1=e^\alpha_{K'} u_2$. Hence, $\rho^{1}_{K,K'} \,q u_1=j^\alpha_{K'} u_2$, where $q\colon \mathcal U(K)\to S^{\prime
0}_1(K)$ is the canonical isomorphism existing because $K$ is proper, i.e., $q u_1$ and $u_2$ have the same restrictions to $S^0_1(K')$. By
Theorem~\ref{t4}, there is a $\hat u\in S^{\prime 0}_\alpha(K)$ such that $u_2$ is the restriction of $\hat u$ to $S^0_\alpha(K')$. Hence $\hat u$ is
an extension of $u$ and, therefore, $u$ is carried by $K$. To prove Theorem~\ref{t4}, we need the next lemma.

\begin{lem}\label{l2}
Let $\alpha\geq 1$ and $K\subset K'$ be proper closed cones in $\R^k$. The space $S^0_1(K')$ is dense in the space $S^0_1(K)\cap S^0_\alpha(K')$
endowed with its natural intersection topology.
\end{lem}
\begin{proof}
Let $g\in S^0_1(K)\cap S^0_\alpha(K')$. Then there are $A,B>0$ and proper conic neighborhoods $W$ and $W'$ of $K$ and $K'$ respectively such that
$g\in S^{0,B}_{1,A}(W)$ and $g\in S^{0,B}_{\alpha,A}(W')$. We can assume that $W\subset W'$ (otherwise we can replace $W$ with $W\cap W'$). Let $f$
be a function on $\C^k$ such that $f(0)=1$ and $f\in S^{0,b}_{1,a}(W')$ for some $a,b>0$ (for example, let $l$ be a linear functional on $\R^k$ such
that $l(x)>0$ for any $x\in \bar W'\setminus\{0\}$; then $f(x+iy)= e^{-l(x)-il(y)}$ belongs to $S^{0,b}_{1,a}(W')$ for $a$ and $b$ large enough). For
$n=1,2,\ldots$, we set $g_n(z)=g(z) f(z/n)$. It follows immediately from Definition~\ref{d1} that
\begin{multline}\label{2:1}
\max\left(\|g_n\|^1_{W,A,\widetilde B},\, \|g_n\|^1_{W',na,\widetilde B},\,\|g_n\|^\alpha_{W',A,\widetilde B}\right) \leq
\\ \leq\|f\|^1_{W',a,b}\,(\|g\|^1_{W,A,B}+\|g\|^\alpha_{W',A,B})
\end{multline}
for all $n$, where $\widetilde B=B+b$. Hence, $g_n\in S^0_1(K')$ and $g_n\in S^{0,\widetilde B}_{1,A}(W)\cap S^{0,\widetilde B}_{\alpha,A}(W')$ for
all $n$. Choose $A'>A$ and $B'>\widetilde B$. To prove the statement, it suffices to show that $g_n\to g$ in both $S^{0,B'}_{1,A'}(W)$ and
$S^{0,B'}_{\alpha,A'}(W')$. By Definition~\ref{d1}, we have
\begin{align}
&|g_n(z)-g(z)|e^{-\sigma^1_{W,A',B'}(z)}\leq \|g_n-g\|^1_{W,A,\widetilde B}\, e^{-\eta|x|-(B'-\widetilde B)|y|},\label{2:2}\\
&|g_n(z)-g(z)|e^{-\sigma^\alpha_{W',A',B'}(z)}\leq \|g_n-g\|^\alpha_{W',A,\widetilde B}\, e^{-\eta|x|^{1/\alpha}-(B'-\widetilde B)|y|},\label{2:3}
\end{align}
where $z=x+iy$ and $\eta$ is the minimum of $A^{-1/\alpha}-A^{\prime -1/\alpha}$ and $A^{-1}-A^{\prime -1}$. For $R>0$, let $Q_R$ be the compact set
$\{x+iy\in\C^k: |x|\leq R,\,|y|\leq R\}$. Fix $\epsilon>0$. By (\ref{2:1}), (\ref{2:2}), and (\ref{2:3}), there exists $R$ such that the left-hand
sides of (\ref{2:2}) and (\ref{2:3}) do not exceed $\epsilon$ for all $n$ and any $z\notin Q_R$. Since $g_n$ converge to $g$ uniformly on compact
sets in $\C^k$, there is $n_0$ such that the the left-hand sides of (\ref{2:2}) and (\ref{2:3}) do not exceed $\epsilon$ for any $n\geq n_0$ and
$z\in Q_R$. Hence $\|g_n-g\|^1_{W,A',B'}\leq\epsilon$ and $\|g_n-g\|^\alpha_{W',A',B'}\leq\epsilon$ for any $n\geq n_0$. The lemma is proved.
\end{proof}

\noindent {\bf Proof of Theorem~\ref{t4}:} Let $l\colon S^0_1(K)\cap S^0_\alpha(K')\to S^0_1(K)\oplus S^0_\alpha(K')$ and $m\colon S^0_1(K) \oplus
S^0_\alpha(K')\to S^0_\alpha(K)$ be the continuous linear mappings defined by the relations $l(f)=(f,-f)$ and $m(f_1,f_2)=f_1+f_2$ (the space
$S^0_1(K)\cap S^0_\alpha(K')$ is endowed with the intersection topology). Clearly, $\mathrm{Im}\,l=\mathrm{Ker}\,m$. Let $v$ be the continuous linear
functional on $S^0_1(K)\oplus S^0_\alpha(K')$ defined by the relation $v(f_1,f_2)=u_1(f_1)+u_2(f_2)$. By the assumption, we have $v(l(f))=0$ for any
$f\in S^0_1(K')$. In view of Lemma~\ref{l2}, this implies that $vl=0$ and hence $\mathrm{Ker}\,m\subset \mathrm{Ker}\,v$. If the mapping $m$ is
surjective, then the open mapping theorem (see Theorem~IV.8.3 in~\cite{Schaefer}; it is applicable because DFS spaces are strong duals of reflexive
Fr\'echet spaces~\cite{Komatsu} and, therefore, are B-complete) implies that $S^0_\alpha(K)$ is topologically isomorphic to the quotient space
$(S^0_1(K)\oplus S^0_\alpha(K'))/\mathrm{Ker}\,m$. It hence follows that there exists a continuous linear functional $u$ on $S^0_\alpha(K)$ such that
$v=u\circ m$. If $f_1\in S^0_1(K)$ and $f_2\in S^0_\alpha(K')$, then we have $u(f_1)=u m(f_1,0)= v(f_1,0)=u_1(f_1)$ and $u(f_2)= u m(0,f_2)=
v(0,f_2)=u_2(f_2)$, i.e., $u_1$ and $u_2$ are the restrictions of $u$ to $S^0_1(K)$ and $S^0_\alpha(K')$ respectively. Proving the statement thus
reduces to proving the surjectivity of $m$. The latter is implied by the following result on the decomposition of test functions.

\begin{thm}\label{t5}
Let $\alpha\geq 1$ and $K\subset K'$ be proper closed cones in $\R^k$. For any $f\in S^0_\alpha(K)$, there exist $f_1\in S^0_1(K)$ and $f_2\in
S^0_\alpha(K')$ such that $f=f_1+f_2$.
\end{thm}

The case $\alpha=1$ is obvious, so we assume $\alpha>1$. The proof of Theorem~\ref{t5} essentially relies on the following H\"ormander's
$L^2$-estimate for the solutions of the inhomogeneous Cauchy--Riemann equations (see Theorem~4.2.6 in~\cite{Hoer}).

\begin{lem}\label{L2est}
Let $\rho$ be a plurisubharmonic function on $\C^k$ and $\eta_j$, $j=1,\ldots,k$, be locally square-integrable functions on $\C^k$. If
\[
\int|\eta_j(z)|^2 e^{-\rho(z)}\,\di\lambda(z)<\infty,\quad j=1,\ldots,k,
\]
where $\di\lambda$ be the Lebesgue measure on $\C^k$, and $\eta_j$ (as generalized functions) satisfy the compatibility conditions\footnote{Here and
hereafter, we use the short notation $\bar\partial_j$ for $\partial/\partial\bar z_j$.} $\bar\partial_j\eta_l=\bar\partial_l\eta_j$, then the
inhomogeneous Cauchy--Riemann equations $\bar\partial_j\psi=\eta_j$ have a locally square-integrable solution satisfying the estimate\footnote{The
estimate in Lemma~\ref{L2est} differs from the estimate in~\cite{Hoer} by the factor $k^2$ in the right-hand side, which appears because we use the
uniform norm instead of the Euclidean norm used in~\cite{Hoer}.}
\[
2\int |\psi(z)|^2e^{-\rho(z)}(1+|z|^2)^{-2}\,\di\lambda(z)
\leq k^2\sum_{j=1}^k \int |\eta_j(z)|^2
e^{-\rho(z)}\,\di\lambda(z).
\]
\end{lem}

The idea of the proof of Theorem~\ref{t5} is as follows. We first construct a \textit{smooth} decomposition $f=\tilde f_1+\tilde f_2$, where $\tilde
f_1$ and $\tilde f_2$ have the same growth properties as elements of $S^0_1(K)$ and $S^0_\alpha(K')$ respectively, and look for $f_1$ and $f_2$ in
the form $f_1=\tilde f_1 -\psi$ and $f_2=\tilde f_2+\psi$. Then the requirement that $f_{1,2}$ be analytic implies the equations $\bar\partial_j
\psi=\bar\partial_j\tilde f_1=-\bar\partial_j\tilde f_2$. We can therefore use Lemma~\ref{L2est} to find their solution that is small enough to
ensure that $f_{1,2}$ have the same growth properties as $\tilde f_{1,2}$ and, hence, satisfy the conditions of Theorem~\ref{t5}.

However, this strategy implies using $L^2$-type norms, while the spaces $S^{0}_{\alpha}(U)$ are defined by supremum norms. To pass to $L^2$-norms, we
make use of results of~\cite{Smirnov}, where this problem was considered for a broad class of spaces containing all spaces $S^0_\alpha(U)$ with
$\alpha\geq 1$. Given $\alpha\geq 1$, $A,B>0$, and a cone $U$ in $\R^k$, let $H^{0,B}_{\alpha,A}(U)$ be the Hilbert space of entire functions on
$\C^k$ with the finite norm
\begin{equation}\label{2:37}
|f|^\alpha_{U,A,B}=\left[\int
|f(z)|^2e^{-2\sigma^\alpha_{U,A,B}(z)}\,\di\lambda(z)\right]^{1/2},
\end{equation}
where  $\sigma^\alpha_{U,A,B}$ is given by~(\ref{000}). It follows from Lemma~4 in~\cite{Smirnov} that
\begin{equation}\label{repr}
S^0_\alpha(U)=\bigcup_{A,B>0,\,W\supset U}H^{0,B}_{\alpha,A}(W),
\end{equation}
where $W$ ranges all conic neighborhoods of $U$ and the union is endowed with the inductive limit topology. The next elementary lemma, which follows
from Lemma~9 in~\cite{Smirnov}, summarizes some simple facts about cones in $\R^k$ needed for the proof of Theorem~\ref{t5}.

\begin{lem}\label{l3a}
Let $K_1$ and $K_2$ be closed cones in $\R^k$ such that $K_1\cap K_2=\{0\}$.
\begin{itemize}

\item[{\rm A.}] There exist conic neighborhoods $V_{1,2}$ of $K_{1,2}$ such that $\bar V_1\cap \bar V_2=\{0\}$.

\item[{\rm B.}] There exists $\theta>0$ such that $\delta_{K_1}(x)\geq \theta|x|$ for any $x\in K_2$.
\end{itemize}
\end{lem}

Given a closed cone $K$ in $\R^k$ and its conic neighborhood $U$, there is a conic neighborhood $V$ of $K$ such that $\bar V\subset U$ (apply
Lemma~\ref{l3a}(A) to $K_1=K$ and $K_2=(\R^k\setminus U)\cup\{0\}$). We shall derive Theorem~\ref{t5} from the next lemma.

\begin{lem}\label{l8}
Let $\alpha>1$ and $A,B>0$. Let $V\subset \R^k$ be a proper cone, $W$ be a proper conic neighborhood of $\bar V$, and $U$ be a proper cone containing
$W$. For any $\epsilon>0$, $A'>A$, and $f\in H^{0,B}_{\alpha,A}(W)$, there exist $B'>0$, $f_1\in H^{0,B'}_{1,\epsilon}(V)$, and $f_2\in
H^{0,B'}_{\alpha,A'}(U)$ such that $f=f_1+f_2$.
\end{lem}

Let $f$ satisfy the conditions of Theorem~\ref{t5} and $U$ be a proper conic neighborhood of $K'$. By~(\ref{repr}), there are $A,B>0$ and a conic
neighborhood $W$ of $K$ such that $f\in H^{0,B}_{\alpha,A}(W)$. We can assume that $W\subset U$ (otherwise we replace $W$ with $U\cap W$). Let $V$ be
a conic neighborhood of $K$ such that $\bar V\subset W$. By Lemma~\ref{l8}, we have $f=f_1+f_2$, where $f_1\in H^{0,B'}_{1,\epsilon}(V)$ and $f_2\in
H^{0,B'}_{\alpha,A'}(U)$ for some $A',B',\epsilon>0$. It now follows from~(\ref{repr}) that $f_1\in S^0_1(K)$ and $f_2\in S^0_\alpha(K')$ and,
therefore, Lemma~\ref{l8} implies Theorem~\ref{t5}.

As explained above, the proof of Lemma~\ref{l8} is based on the $L^2$-estimate given by Lemma~\ref{L2est}. However, Lemma~\ref{L2est} involves
plurisubharmonic functions, while the weight functions $\sigma^\alpha_{U,A,B}$ used in the definition of the spaces $H^{0,B}_{\alpha,A}(U)$ are not
plurisubharmonic. This problem is resolved by the next lemma.

\begin{lem}\label{l7}
Let $\alpha>1$, $U$ be a proper cone in $\R^k$, and $K_1$ and $K_2$ be closed cones such that $K_1\cap K_2=\{0\}$ and $K_1\cup K_2\subset U$. For any
$\kappa,d\geq 0$, there exist constants $b>0$ and $H$ and a plurisubharmonic function $\varrho$ on $\C^k$ such that
\begin{align}
&\varrho(x+iy)\leq -|x|^{1/\alpha}+b\,\delta_U(x)+b|y|,\quad x,y\in\R^k,\label{2:30}\\
&\varrho(x+iy)\leq -\kappa|x|+b\,\delta_{K_1}(x)+b|y|,\quad x,y\in\R^k,\label{2:31}\\
&\varrho(x+iy)\geq -|x|^{1/\alpha}-H,\quad x\in K_2^d,\, y\in\R^k, \label{2:32}
\end{align}
where $K_2^d=\{x\in\R^k: \delta_{K_2}(x)\leq d\}$ is the closed $d$-neighborhood of $K_2$.
\end{lem}

Before proving Lemma~\ref{l7}, we derive Lemma~\ref{l8} from Lemma~\ref{l7}.

\begin{proof}[Proof of Lemma~$\ref{l8}$]
Without loss of generality, we can assume that $\bar W\subset U$ (otherwise we can replace $U$ with $\bar U$). Fix $\delta>0$ and choose a
nonnegative smooth function $g_0$ on $\R^k$ such that $g_0(x) = 0$ for $|x|\geq \delta$ and $\int_{\R^k} g_0(x)\,\di x=1$. We define the smooth
functions $g_1$ and $g_2$ on $\R^k$ by the relations
\[
g_1(x)=\int_{\R^k\setminus W}g_0(x-\xi)\,\di\xi,\quad g_2(x)=\int_W g_0(x-\xi)\,\di\xi,\quad x\in\R^k.
\]
For any $x\in\supp g_2$, we have $\delta_W(x)\leq \delta$. Hence,
\begin{equation}\label{2:38}
\sigma^\alpha_{W,A,B}(z)\leq \sigma^\alpha_{\R^k,A,B}(z) + B\delta,\quad z\in\supp g_2+i\R^k.
\end{equation}
By Lemma~\ref{l3a}(A), there exists a conic neighborhood $W'$ of $(\R^k\setminus W)\cup\{0\}$ such that $\bar W'\cap \bar V=\{0\}$. By
Lemma~\ref{l3a}(B), there is $\theta>0$ such that $\delta_{\R^k\setminus W}(x)\geq \theta|x|$ for any $x\notin W'$. For $x\in \supp g_1$, we have
$\delta_{\R^k\setminus W}(x)\leq \delta$ and, therefore, $\supp g_1\setminus W'$ is a bounded set. It follows from Lemma~\ref{l3a}(B) that
$\delta_V(x)\geq \theta'|x|$ for some $0<\theta'\leq 1$ and any $x\in W'$. Since $\delta_W(x)\leq |x|$, we have
\begin{equation}\label{2:39}
\sigma^\alpha_{W,A,B}(z)\leq \sigma^1_{V,\epsilon,\tilde B}(z) + R,\quad z\in\supp g_1+i\R^k,
\end{equation}
where $R$ is a constant and $\tilde B = B/\theta'+1/(\theta'\epsilon)$. We define the smooth functions $\tilde f_1$ and $\tilde f_2$ on $\C^k$ by the
relation $\tilde f_{1,2}(x+iy)=g_{1,2}(x) f(x+iy)$. Since $f$ is analytic and $g_1+g_2=1$, we have
\begin{equation}\label{2:40}
\bar\partial_j \tilde f_1(z)=-\bar\partial_j \tilde f_2(z)=
-\frac{1}{2}f(z)\frac{\partial g_2(x)}{\partial x_j},\quad j=1,\ldots,k,\,\,z=x+iy\in\C^k.
\end{equation}
It follows from the definition of $g_2$ that all its partial derivatives are bounded on $\R^k$, and in view of~(\ref{2:38}), (\ref{2:39}),
and~(\ref{2:40}), we obtain
\begin{equation}\label{2:41}
|\tilde f_1|^1_{V,\epsilon,\tilde B}<\infty,\quad |\tilde f_2|^\alpha_{\R^k,A,B}<\infty, \quad |\bar\partial_j\tilde f_1|^\alpha_{\R^k,A,B}<\infty
\end{equation}
for any $j=1,\ldots,k$. We now choose $\kappa>A/\epsilon$ and set $d=\delta/A$, $K_1=\bar V$ and $K_2=\partial W$ (the boundary of $W$). By
Lemma~\ref{l7}, there is a plurisubharmonic function $\varrho$ on $\C^k$ satisfying (\ref{2:30})--(\ref{2:32}). Let
\begin{equation}\label{2:42}
\sigma(z)=\varrho(z/A)+B|y|+H,\quad z=x+iy\in\C^k,
\end{equation}
where $H$ is the constant entering~(\ref{2:32}). Clearly, $\sigma$ is a plurisubharmonic function on $\C^k$, and it follows from (\ref{2:32}) and
(\ref{2:42}) that $\sigma(x+iy)\geq \sigma^\alpha_{\R^k,A,B}(x+iy)$ for any $x\in K^\delta_2$ and $y\in\R^k$. Because $g_1+g_2=1$, we have
$\supp\partial_j g_2\subset \supp g_1\cap\supp g_2\subset K^\delta_2$, and in view of (\ref{2:40}) and (\ref{2:41}), we obtain
\begin{equation}\label{2:43}
\int |\bar\partial_j\tilde f_1(z)|^2 e^{-2\sigma(z)}\,\di\lambda(z)\leq (|\bar\partial_j\tilde f_1|^\alpha_{\R^k,A,B})^2<\infty,\quad j=1,\ldots,k.
\end{equation}
Let $\tilde\sigma(z)=\sigma(z)+\log(1+|z|^2)$. By Lemma~\ref{L2est}, the equations $\bar\partial_j\psi=\bar\partial_j\tilde f_1$ have a locally
square-integrable solution such that
\begin{equation}\label{2:44}
\int |\psi(z)|^2 e^{-2\tilde\sigma(z)}\,\di\lambda(z) < \infty.
\end{equation}
Let $\tilde b> B+b/A$. It follows from~(\ref{2:30}), (\ref{2:31}), and~(\ref{2:42}) that
\begin{equation}\label{2:45}
\tilde\sigma(z)\leq \sigma^1_{V,\epsilon,\tilde b}(z)+C, \qquad \tilde\sigma(z)\leq \sigma^\alpha_{U,A',\tilde b}(z)+C, \qquad z\in\C^k,
\end{equation}
where $C$ is a constant. In view of~(\ref{2:40}), we have $\bar\partial_j(\tilde f_1-\psi)=\bar\partial_j(\tilde f_2+\psi)=0$; hence, there are
entire functions $f_1$ and $f_2$ that coincide almost everywhere with $\tilde f_1-\psi$ and $\tilde f_2+\psi$ respectively. By~(\ref{2:44})
and~(\ref{2:45}), we have $|\psi|^1_{V,\epsilon,\tilde b}<\infty$ and $|\psi|^\alpha_{U,A',\tilde b}<\infty$. In view of~(\ref{2:41}), this implies
that $f_1\in H^{0,B'}_{1,\epsilon}(V)$ and $f_2\in H^{0,B'}_{\alpha,A'}(U)$ for any $B'\geq\max(\tilde B,\tilde b)$. Moreover, $f = f_1+f_2$ because
continuous functions coinciding almost everywhere are equal. The lemma is proved.
\end{proof}

The rest of this section is devoted to the proof of Lemma~\ref{l7}. Let $\Theta$ be the subharmonic function defined by the relation
\begin{equation}\label{theta}
\Theta(z)=\log\left|\frac{\sin z}{z}\right|,\quad z\in \C.
\end{equation}

\begin{lem}\label{l4}
The function $\Theta$ is strictly decreasing on the segment $[0,\pi]$ of the real axis and satisfies the inequalities
\begin{align}
&\Theta(x+iy)\geq \Theta(x),\quad x,y\in\R. \label{2:12} \\
&\Theta(x+iy)\leq \Theta(x) + |y|,\quad x,y\in\R,\,\,|x|\leq 3\pi/4. \label{2:13}
\end{align}
\end{lem}
\begin{proof}
We have $(\sin x/x)'=x^{-2}(x\cos x -\sin x)$. The function $x\cos x -\sin x$ vanishes at $x=0$ and strictly decreases on $[0,\pi]$ because its
derivative $-x\sin x$ is strictly negative for $0<x<\pi$. This implies that
\begin{equation}\label{ineq}
x\cos x < \sin x, \quad 0< x\leq\pi,
\end{equation}
and, therefore, $\sin x/x$ strictly decreases on $[0,\pi]$. Hence, $\Theta$ strictly decreases on $[0,\pi]$. It is straightforward to check that
\begin{equation}\label{2:14}
\left|\frac{\sin z}{z}\right|^2 = \frac{\sin^2 x + \sinh^2 y}{x^2+y^2},\quad z=x+iy\in \C.
\end{equation}
Let $0\leq s< t$. Since the function $(u+s)(u+t)^{-1}$ is increasing in $u$ on $[0,\infty)$, we have $(u+s)(u+t)^{-1}\geq s/t$ for any $u\geq 0$.
Setting $s=\sin^2 x$, $t=x^2$, and $u=\sinh^2 y$ and applying this inequality, we derive from~(\ref{2:14}) that
\[
\left|\frac{\sin z}{z}\right|^2\geq \frac{\sin^2 x + \sinh^2 y}{x^2+\sinh^2 y}\geq \frac{\sin^2 x}{x^2}, \quad z=x+iy\in\C,\, x\ne 0.
\]
By continuity, this inequality remains valid for $x=0$, and passing to the logarithms, we obtain~(\ref{2:12}). Further, it easily follows
from~(\ref{2:14}) that
\begin{equation}\nonumber
\Theta(x+iy)-\Theta(x) = |y| +\frac{1}{2}\log \frac{\sin^2 x+\cos 2x \frac{1-e^{-2|y|}}{2}-\frac{1-e^{-4|y|}}{4}}{\sin^2 x+\frac{\sin^2 x}{x^2} y^2}.
\end{equation}
Hence, to prove~(\ref{2:13}), it suffices to show that
\begin{equation}\label{2:15}
\cos 2x \frac{1-e^{-2|y|}}{2}-\frac{1-e^{-4|y|}}{4} \leq \frac{\sin^2 x}{x^2} y^2,\quad \quad x,y\in\R,\,\,|x|\leq 3\pi/4.
\end{equation}
If $\pi/4\leq |x|\leq 3\pi/4$, then $\cos 2x \leq 0$ and (\ref{2:15}) is obvious. For any $x,y\in\R$, we have $1-e^{-4|y|}\geq \cos 2x (1-e^{-4|y|})$
and, therefore,
\[
\cos 2x \frac{1-e^{-2|y|}}{2}-\frac{1-e^{-4|y|}}{4}\leq \cos 2x
\left(\frac{1-e^{-2|y|}}{2}\right)^2,\quad x,y\in\R.
\]
Since $(1-e^{-2|y|})/2\leq |y|$ for any $y\in\R$ and $\cos 2x\leq \cos x$ for $|x|\leq \pi/2$, inequality~(\ref{2:15}) will be proved if we
demonstrate that $\cos x\leq x^{-2}\sin^2 x$ for $|x|\leq \pi/4$. This latter inequality holds because the function $x^2\cos x - \sin^2 x$ vanishes
at $x=0$ and decreases on $[0,\pi]$ (in view of~(\ref{ineq}), its derivative does not exceed $-4\sin x((x/2)^2-\sin^2(x/2))$ and, hence, is
nonpositive for $0\leq x\leq \pi$). The lemma is proved.
\end{proof}

We define the function $\mu$ on $[0,\infty)$ by the relation
\[
\mu(x)=\left\{
\begin{matrix}
\Theta(x),& 0\leq x\leq \pi/2,\\
-\log|x|,& x>\pi/2.
\end{matrix}
\right.
\]
Thus, $\mu$ is a continuous function on $[0,\infty)$ such that $\mu(0)=0$ and $\mu(x)<0$ for $x\ne 0$. It follows from Lemma~\ref{l4} that $\mu$
strictly decreases on $[0,\infty)$ and
\begin{equation}\label{2:16}
\Theta(x+iy)\geq \mu(|x|),\quad x,y\in\R,\,\,|x|\leq \pi/2.
\end{equation}
For any $x,y\in\R$, we have $|\sin(x+iy)|\leq e^{|y|}$. Hence, $\Theta(x+iy)\leq -\log|x|+|y|$ and in view of Lemma~\ref{l4}, we obtain
\begin{equation}\label{2:17}
\Theta(x+iy)\leq \mu(|x|)+|y|,\quad x,y\in\R.
\end{equation}

\begin{lem}\label{l5}
Let $0<a<b$ and let $\chi$ be the characteristic function of the segment $[a,b]$ (i.e., $\chi(x)=1$ for $x\in [a,b]$ and $\chi(x)=0$ for $x\notin
[a,b]$). For any $\varkappa>0$, there exist $R>0$, $x_0\in[a,b]$, and a subharmonic function $\Psi$ on $\C$ such that
\begin{equation}\label{2:18}
\begin{aligned}
&\Psi(x+iy)\leq x\chi(x)+R|y|,\\
&\Psi(x_0+iy)\geq x_0,
\end{aligned}
\qquad x,y\in \R,
\end{equation}
and $\Psi$ is bounded below on the strip $\{x+iy\in\C: |x|\leq \varkappa\}$.
\end{lem}
\begin{proof}
For $z\in\C$, let
\[
t(z)=\frac{\pi}{2(b+\varkappa)}\left(z-\frac{b+a}{2}\right).
\]
We define the function $\tilde \mu$ on $\R$ by the relation $\tilde\mu(x)=\mu(|t(x)|)-\mu(h)$, where $h=t(b)=-t(a)$, and set
$\lambda=\sup_{x\in[a,b]} \tilde\mu(x)/x$. The function $\tilde\mu$ is continuous, and $\tilde\mu(x)>0$ for $a<x<b$. Hence $\lambda>0$ and there
exists $x_0\in [a,b]$ such that
\begin{equation}\label{2:19}
\lambda=\tilde\mu(x_0)/x_0.
\end{equation}
Since $\tilde\mu(x)\leq 0$ for $x\notin (a,b)$, we have
\begin{equation}\label{2:20}
\tilde\mu(x)\leq \lambda x\,\chi(x),\quad x\in\R.
\end{equation}
We now set $\Psi(z)=\lambda^{-1}[\Theta(t(z))-\mu(h)]$. Substituting $t(x+iy)$ for $x+iy$ in~(\ref{2:17}), we obtain $\Psi(z)\leq
\lambda^{-1}\tilde\mu(x)+R|y|$ for any $x,y\in\R$, where $R=\pi/[2\lambda(b+\varkappa)]$. In view of~(\ref{2:20}), this implies the upper inequality
in~(\ref{2:18}). Since $|t(x_0)|\leq h\leq \pi/2$, substituting $t(x_0+iy)$ for $x+iy$ in (\ref{2:16}) yields $\Psi(x_0+iy)\geq
\lambda^{-1}\tilde\mu(x_0)$ for any $y\in\R$. Together with~(\ref{2:19}), this gives the lower inequality in~(\ref{2:18}). If $|x|\leq \varkappa$,
then $|t(x)|\leq \pi/2$. By~(\ref{2:16}) and the monotonicity of $\mu$, it hence follows that $\Psi(x+iy)\geq \lambda^{-1}(\mu(\pi/2)-\mu(h))$ for
any $x,y$ such that $|x|\leq \varkappa$. The lemma is proved.
\end{proof}

\begin{lem}\label{l6}
Let $K\subset \R^k$ be a proper closed cone, $V$ be its conic neighborhood, and $l$ be a linear functional on $\R^k$ such that $K\setminus\{0\}$ is
contained in the open halfspace $\{x\in\R^k: l(x)>0\}$. Then there exist constants $r,r'\geq 0$ and a plurisubharmonic function $\Phi$ on $\C^k$ such
that
\begin{align}
&-r'|x|\leq\Phi(x+iy)\leq \max(l(x),0)+r|y|,\quad x,y\in\R^k,\label{2:21}\\
&\Phi(x+iy)\leq r|y|,\quad x\notin \bar V, \,y\in\R^k,\label{2:22}\\
&\Phi(x+iy)\geq l(x),\quad x\in K,\, y\in\R^k. \label{2:23}
\end{align}
\end{lem}
\begin{proof}
Without loss of generality, we can assume $K\ne\{0\}$ (otherwise we can set $\Phi=0$). If $k=1$, then we have either $K=\bar\R_+$ or $K=\bar\R_-$,
and $\Phi(x+iy)=\max(l(x),0)$ satisfies (\ref{2:21})--(\ref{2:23}) with $r=r'=0$. From now on, we assume $k>1$. Let $\lambda=\inf_{x\in
K,\,|x|=1}l(x)$. Since the infimum is taken over a compact set, where $l$ is strictly positive, we have $\lambda>0$. We thus obtain
\begin{equation}\label{2:24}
|x|\leq l(x)/\lambda,\quad x\in K.
\end{equation}
Let $Q=K\cap \{\xi\in\R^k: l(\xi)=1\}$. By~(\ref{2:24}), $Q$ is bounded and, therefore, compact. Choose a basis $e_1,\ldots, e_{k-1}$ in
$\mathrm{Ker}\,l$. For $\xi\notin \mathrm{Ker}\,l$, let $l^1_\xi,\ldots,l^{k-1}_\xi$ be linear functionals on $\R^k$ such that $l^i_\xi(\xi)=0$ and
$l^i_\xi(e_j)=\delta^i_j$ (in other words, $l,l^1_\xi,\ldots,l^{k-1}_\xi$ is the dual basis of $\xi/l(\xi),e_1,\ldots,e_{k-1}$). For any $\xi\notin
\mathrm{Ker}\,l$, we define the norm $|\cdot|_\xi$ on $\R^k$ by the relation
\[
|x|_\xi=|l(x)|+|l^1_\xi(x)|+\ldots+|l^{k-1}_\xi(x)|,\quad x\in\R^k.
\]
Let $M$ and $m$ be, respectively, the supremum and infimum of $|x|_\xi$ on the compact set $\{(x,\xi)\in \R^{2k}: |x|=1,\,\xi\in Q\}$. Since
$|x|_\xi$ is strictly positive and continuous\footnote{The continuity follows from the fact that the mapping $\xi\to l^i_\xi$ from
$\R^k\setminus\mathrm{Ker}\,l$ to the space of linear functionals on $\R^k$ is continuous for any $i=1,\ldots,k-1$.} on this set, we conclude that
$0<m\leq M<\infty$. We therefore have
\begin{equation}\label{2:25}
m|x|\leq |x|_\xi\leq M|x|,\quad \xi\in Q,\,x\in\R^k.
\end{equation}
Let $0<a<b$ and $\chi$ be the characteristic function of $[a,b]$. By Lemma~\ref{l5}, there are $R>0$, $x_0\in [a,b]$, and a subharmonic function
$\Psi$ on $\C$ such that inequalities~(\ref{2:18}) hold and $\Psi$ is bounded below on the strip $|x|\leq 1$. Given a linear functional $L$ on
$\R^k$, we denote by $\hat L$ its unique complex-linear extension to $\C^k$: $\hat L(x+iy)=L(x)+iL(y)$. For $\xi\in Q$ and $\tau>0$, we set
\begin{equation}\label{2:26}
\Phi_{\xi,\tau}(z)=\Psi(\hat l(z))+\tau\sum_{j=1}^{k-1}\Theta\left(\hat l^j_\xi(z)\right),\quad z\in\C^k,
\end{equation}
where $\Theta$ is given by~(\ref{theta}). Further, we set
\begin{equation}\label{2:26a}
\widetilde \Phi_\tau(z) = \sup_{\xi\in Q,\,s>0} s\,\Phi_{\xi,\tau}(z/s),\quad \Phi_\tau(z)=\varlimsup_{z'\to z} \widetilde \Phi_\tau(z').
\end{equation}
Clearly, $\Phi_{\xi,\tau}$ is a plurisubharmonic function on $\C^k$ for any $\xi\in Q$ and $\tau>0$. Hence $\Phi_\tau$ is also a plurisubharmonic
function (see Sec.~II.10.3 in~\cite{V}). We shall show that $\Phi=\Phi_\tau$ satisfies (\ref{2:21}), (\ref{2:22}), and (\ref{2:23}) for some
$r,r'\geq0$ if $\tau$ is large enough. In view of~(\ref{2:17}) and~(\ref{2:18}), it follows from (\ref{2:26}) that
\begin{multline}\nonumber
\Phi_{\xi,\tau}(x+iy)\leq l(x)\chi(l(x))+R|l(y)|+\tau \sum_{j=1}^{k-1}\left(|l^j_\xi(y)|+\mu(|l^j_\xi(x)|)\right)\leq\\
\leq l(x)\chi(l(x))+(R+\tau)|y|_\xi+ \tau \mu\left(\frac{|l^1_\xi(x)|+\ldots+|l^{k-1}_\xi(x)|}{k-1}\right),\quad x,y\in\R^k.
\end{multline}
where we have used the monotonicity and nonpositivity of $\mu$. Note that $l(x-l(x)\xi)=0$ and $l^i_\xi(x-l(x)\xi)=l^i_\xi(x)$ for any $x\in \R^k$.
This implies that $|l^1_\xi(x)|+\ldots+|l^{k-1}_\xi(x)|=|x-l(x)\xi|_\xi$, and using~(\ref{2:25}), we obtain
\begin{equation}\label{2:27}
\Phi_{\xi,\tau}(x+iy)\leq l(x)\chi(l(x)) + r_\tau|y|+\tau\mu\left(\frac{m}{k-1}|x-l(x)\xi|\right),\quad x,y\in\R^k,
\end{equation}
where $r_\tau=M(R+\tau)$. Since $x\chi(x)\leq \max(x,0)$ for any $x\in\R$ and $\mu$ is nonpositive, it follows from (\ref{2:27}) that
$\Phi_{\xi,\tau}(x+iy)\leq \max(l(x),0) + r_\tau|y|$ for any $\xi\in Q$. This implies that $\Phi_\tau$ satisfies the right inequality in (\ref{2:21})
for $r=r_\tau$. Let $H$ be such that $\Psi(x+iy)\geq -H$ for $|x|\leq 1$. By~(\ref{2:16}) and~(\ref{2:26}), we have $\Phi_{\xi,\tau}(x+iy)\geq
-h_\tau$, where $h_\tau=H-\tau(k-1)\mu(1)$, for any $\xi\in Q$ and $x,y\in\R^k$ such that $|x|_\xi\leq 1$. Let $x\ne 0$ and $s=|x|_\xi$ for some
$\xi\in Q$. By~(\ref{2:25}), we obtain $s\Phi_{\xi,\tau}((x+iy)/s)\geq -Mh_\tau |x|$ for any $y\in\R^k$. In view of (\ref{2:26a}), this ensures the
left inequality in (\ref{2:21}) for $\Phi=\Phi_\tau$ and $r'=Mh_\tau$. Further, it follows immediately from~(\ref{2:27}) that
\begin{equation}\label{2:28}
\Phi_{\xi,\tau}(x+iy)\leq r_\tau|y|,\quad x,y\in \R^k,\,l(x)\notin [a,b],
\end{equation}
for any $\xi\in Q$. Let $S=K\cap\{x\in\R^k:a\leq l(x)\leq b\}$. By~(\ref{2:24}), $S$ is a compact set and, therefore, the distance $d$ between $S$
and the closed set $(\R^k\setminus V)\cup\{0\}$ is strictly positive. If $a\leq l(x)\leq b$, then $l(x)\xi\in S$ for any $\xi\in Q$, and~(\ref{2:27})
yields
\begin{equation}\label{2:29}
\Phi_{\xi,\tau}(x+iy)\leq r_\tau|y|+b+\tau\mu\left(\frac{md}{k-1}\right),\quad y\in\R^k, x\notin V,\,a\leq l(x)\leq b.
\end{equation}
Together with~(\ref{2:28}), this implies that $\Phi_{\xi,\tau}(x+iy)\leq r_\tau|y|$ for any $\xi\in Q$, $x\notin V$, and $y\in \R^k$ if $\tau\geq
-b\,\mu(md/(k-1))^{-1}$. Since $\R^k\setminus\bar V$ is an open set, it now follows from~(\ref{2:26a}) that $\Phi_\tau$ satisfies (\ref{2:22}) for
$\tau$ large enough and $r=r_\tau$. It remains to prove (\ref{2:23}). For $x\in K\setminus\{0\}$, we set $\xi_x=x/l(x)$ and $s_x=l(x)/x_0$. Then we
have $l^i_{\xi_x}(x)=0$ for $i=1,\ldots,k-1$. By~(\ref{2:16}), $\Theta(iy)\geq 0$ for any $y\in\R$, and it follows from~(\ref{2:18}) and (\ref{2:26})
that
\[
s_x \Phi_{\xi_x,\tau}((x+iy)/s_x)\geq s_x \Psi(x_0+il(y)/s_x)\geq s_x x_0=l(x),\quad x\in K\setminus\{0\},\,y\in\R^k.
\]
In view of (\ref{2:26a}), this implies that $\Phi_\tau$ satisfies (\ref{2:23}). The lemma is proved.
\end{proof}

\begin{proof}[Proof of Lemma~$\ref{l7}$]
Since every proper cone has a proper conic neighborhood, we can assume that $U\setminus\{0\}$ is an open set. By Lemma~\ref{l3a}(A), there exist
conic neighborhoods $V_1$ and $V_2$ of $K_1$ and $K_2$ respectively such that $V_1\cup V_2\subset U$ and $\bar V_1\cap \bar V_2 = \{0\}$. Let $l$ be
a linear functional on $\R^k$ such that $l(x)>0$ for any $x\in\bar U\setminus\{0\}$ and
\begin{equation}\label{2:32a}
\inf_{x\in \bar V_1,\,|x|=1} l(x)\geq\kappa.
\end{equation}
By Lemma~\ref{l6}, there exist $r,r'\geq 0$ and a plurisubharmonic function $\Phi$ such that (\ref{2:21})--(\ref{2:23}) hold for $K=\bar V_2$ and
$V=(\R^k\setminus \bar V_1)\cup\{0\}$. Since the space $S^0_\alpha(\R)$ is nontrivial, Lemma~5 in~\cite{Smirnov} ensures that there are constants
$B>0$ and $H$ and a plurisubharmonic function $\sigma$ on $\C^k$ such that
\begin{equation}\label{2:33}
-|x|^{1/\alpha}-H \leq \sigma(x+iy)\leq -|x|^{1/\alpha}+B|y|,\quad x,y\in\R^k.
\end{equation}
We now define the plurisubharmonic function $\varrho$ by the relation
\begin{equation}\label{2:34}
\varrho(z) = \Phi(z)+\sigma(z)-l(x),\quad z=x+iy\in\C^k.
\end{equation}
By~(\ref{2:23}) and (\ref{2:33}), $\varrho$ satisfies~(\ref{2:32}) for $x\in \bar V_2$. By Lemma~\ref{l3a}(B), we have $\delta_{K_2}(x)\geq
\theta|x|$ for some $\theta>0$ and any $x\notin V_2$ and, hence, $K^d_2\setminus \bar V_2$ is a bounded set. In view of~(\ref{2:21}), (\ref{2:33}),
and (\ref{2:34}), it follows that $\varrho$ is bounded below on the set $\{x+iy\in\C^k: x\in K^d_2\setminus \bar V_2\}$. We can hence
ensure~(\ref{2:32}) for all $x\in K^d_2$ increasing, if necessary, the constant $H$. By~(\ref{2:22}), (\ref{2:32a}), and (\ref{2:33}), we have
\begin{equation}\label{2:35}
\varrho(x+iy)\leq -\kappa|x|+(r+B)|y|,\quad x\in V_1,\,y\in\R^k.
\end{equation}
By Lemma~\ref{l3a}(B), there is $\theta'>0$ such that $\theta'|x|\leq \delta_{K_1}(x)$ for any $x\notin V_1$. It follows from~(\ref{2:21}),
(\ref{2:33}), and (\ref{2:34}) that
\[
\varrho(x+iy)\leq -\kappa|x|+\theta^{\prime -1}(\kappa+|l|)\,\delta_{K_1}(x)+(r+B)|y|,\quad x\notin V_1,\,y\in\R^k,
\]
where $|l|=\sup_{|x|=1}|l(x)|$. Together with~(\ref{2:35}), this estimate implies~(\ref{2:31}) for any $b\geq r+B+(\kappa+|l|)/\theta'$. Further, it
follows from~(\ref{2:21}), (\ref{2:33}), and (\ref{2:34}) that
\begin{equation}\label{2:36}
\varrho(x+iy)\leq \max(-l(x),0)-|x|^{1/\alpha}+(r+B)|y|,\quad x,y\in\R^k.
\end{equation}
Using Lemma~\ref{l3a}(B), it is easy to show that $\max(-l(x),0)\leq b\,\delta_U(x)$ for any $x\in\R^k$ and some $b>0$. Hence (\ref{2:36}) implies
(\ref{2:30}) for $b$ large enough. The lemma is proved.
\end{proof}

\section{Proof of Theorem~\ref{t1}}
\label{s3}

In~\cite{AIM}, the proof of Theorem~\ref{t3} fell into two largely independent parts: the analytic part concerning proper cones and the algebraic
part concerning passing from proper cones to the general case. Here, the situation is much the same, and the problem of deriving Theorem~\ref{t1}
from Theorem~\ref{t4} can be reformulated in a purely algebraic way in terms of abstract inductive systems indexed by partially ordered sets of
certain type.

We first recall some notation and definitions introduced in~\cite{AIM}. By an inductive system $\mathcal X$ of vector spaces indexed by a partially
ordered set $\Gamma$, we mean a family $\{\mathcal X(\gamma)\}_{\gamma\in \Gamma}$ of vector spaces together with a family of linear mappings
$\rho^{\mathcal X}_{\gamma\gamma'}\colon\mathcal X(\gamma)\to \mathcal X(\gamma')$ defined for $\gamma\leq \gamma'$ and such that $\rho^{\mathcal
X}_{\gamma\gamma}$ is the identity mapping for any $\gamma\in \Gamma$ and $\rho^{\mathcal X}_{\gamma\gamma''}=\rho^{\mathcal
X}_{\gamma'\gamma''}\rho^{\mathcal X}_{\gamma\gamma'}$ for $\gamma\leq\gamma'\leq\gamma''$. Let $\iota^{\mathcal X}_\gamma$ denote the canonical
embedding of $\mathcal X(\gamma)$ in $\oplus_{\gamma'\in \Gamma} \mathcal X(\gamma')$. The inductive limit $\varinjlim\mathcal X$ is by definition
the quotient space $[\oplus_{\gamma\in \Gamma} \mathcal X(\gamma)]/N^{\mathcal X}$, where $N^{\mathcal X}$ is the subspace of $\oplus_{\gamma\in
\Gamma} \mathcal X(\gamma)$ spanned by all elements of the form $\iota^{\mathcal X}_\gamma x- \iota^{\mathcal X}_{\gamma'}\rho^{\mathcal
X}_{\gamma\gamma'}x$, $x\in \mathcal X(\gamma)$. The canonical mapping $\rho^{\mathcal X}_\gamma\colon \mathcal X(\gamma)\to \varinjlim\mathcal X$ is
defined by the relation $\rho^{\mathcal X}_\gamma=j^\mathcal X \iota^{\mathcal X}_\gamma$, where $j^\mathcal X$ is the canonical surjection of
$\oplus_{\gamma\in \Gamma} \mathcal X(\gamma)$ onto $\varinjlim\mathcal X$. As in~\cite{AIM}, we do not assume that the index set $\Gamma$ is
directed. It is important that the standard universal property of inductive limits remains valid for such generalized inductive systems.

Recall that a partially ordered set $\Gamma$ is called a lattice if each two-element subset $\{\gamma_1,\gamma_2\}$ of $\Gamma$ has a supremum
$\gamma_1\vee \gamma_2$ and an infimum $\gamma_1\wedge \gamma_2$. A lattice $\Gamma$ is called distributive if $\gamma_1\wedge(\gamma_2\vee\gamma_3)=
(\gamma_1\wedge\gamma_2)\vee(\gamma_1\wedge\gamma_3)$ for any $\gamma_1,\gamma_2,\gamma_3\in \Gamma$.

\begin{defn} \label{dxx4}
A partially ordered set $\Gamma$ is called a quasi-lattice if every two-element subset of $\Gamma$ has an infimum and every bounded above two-element
subset of $\Gamma$ has a supremum. A quasi-lattice $\Gamma$ is called distributive if $\gamma_1\wedge(\gamma_2\vee\gamma_3)=
(\gamma_1\wedge\gamma_2)\vee(\gamma_1\wedge\gamma_3)$ for every bounded above pair $\gamma_2,\gamma_3\in \Gamma$ and every $\gamma_1\in \Gamma$.
\end{defn}

Clearly, every (distributive) lattice is a (distributive) quasi-lattice.

\begin{defn}\label{dxx6}
An inductive system $\mathcal X$ of vector spaces indexed by a quasi-lattice $\Gamma$ is called prelocalizable if the following conditions are
satisfied:
\begin{enumerate}
\item[(I)] The mappings $\rho^{\mathcal X}_{\gamma\gamma'}$ are injective for any $\gamma,\gamma'\in \Gamma$, $\gamma\leq\gamma'$.

\item[(II)] If a pair $\gamma_1,\gamma_2\in \Gamma$ is bounded above and $x\in \mathcal X(\gamma_1\vee\gamma_2)$, then there are $x_{1,\,2}\in
\mathcal X(\gamma_{1,\,2})$ such that $x=\rho^{\mathcal X}_{\gamma_1,\,\gamma_1\vee\gamma_2}(x_1)+ \rho^{\mathcal
X}_{\gamma_2,\,\gamma_1\vee\gamma_2}(x_2)$.

\item[(III)] If a pair $\gamma_1,\gamma_2\in \Gamma$ is bounded above by an element $\gamma\in \Gamma$, $x_{1,\,2}\in\mathcal X(\gamma_{1,\,2})$, and
$\rho^{\mathcal X}_{\gamma_1,\,\gamma}(x_1)= \rho^{\mathcal X}_{\gamma_2,\,\gamma}(x_2)$, then there is an $x\in \mathcal X(\gamma_1\wedge\gamma_2)$
such that $x_1=\rho^{\mathcal X}_{\gamma_1\wedge\gamma_2,\,\gamma_1}(x)$ and $x_2=\rho^{\mathcal X}_{\gamma_1\wedge\gamma_2,\,\gamma_2}(x)$.
\end{enumerate}
\end{defn}

Let $\mathcal X$ be an inductive system indexed by $\Gamma$. For $I\subset \Gamma$, we define the inductive system $\mathcal X^I$ over $I$ setting
$\mathcal X^I(\gamma)=\mathcal X(\gamma)$ and $\rho^{\mathcal X^I}_{\gamma\gamma'}=\rho^{\mathcal X}_{\gamma\gamma'}$ for $\gamma,\gamma'\in I$,
$\gamma\leq \gamma'$ (i.e., $\mathcal X^I$ is the ``restriction'' of $\mathcal X$ to $I$). Let $I\subset J\subset \Gamma$. By the universal property
of inductive limits, $\rho_\gamma^{\mathcal X^J}$ uniquely determine a map $\tau^{\mathcal X}_{I,\,J}\colon \varinjlim \mathcal X^I\to \varinjlim
\mathcal X^J$ such that $\tau^{\mathcal X}_{I,\,J}\rho_\gamma^{\mathcal X^I}= \rho_\gamma^{\mathcal X^J}$ for any $\gamma\in I$. Let $\lambda$ be a
nondecreasing map from $\Gamma$ to a partially ordered set $\Delta$. With every $\delta \in \Delta$, we associate the set $\Gamma_\delta=\{\gamma\in
\Gamma\,|\,\lambda(\gamma)\leq \delta\}$ and define the inductive system $\lambda(\mathcal X)$ over $\Delta$ setting $\lambda(\mathcal
X)(\delta)=\varinjlim \mathcal X^{\Gamma_\delta}$ and $\rho^{\lambda(\mathcal X)}_{\delta\delta'}= \tau^{\mathcal
X}_{\Gamma_\delta,\,\Gamma_{\delta'}}$ for $\delta,\delta'\in \Delta$, $\delta\leq\delta'$.

Let $\mathcal K(\R^k)$ denote the set of all nonempty closed cones in $\R^k$ ordered by inclusion. Clearly, $\mathcal K(\R^k)$ is a distributive
lattice, while the set $\mathcal P(\R^k)$ of closed proper cones in $\R^k$ is a distributive quasi-lattice. For any $\alpha\geq 1$, the spaces
$S^{\prime 0}_\alpha(K)$, $K\in \mathcal K(\R^k)$, together with the canonical mappings $\rho^\alpha_{K,K'}\colon S^{\prime 0}_\alpha(K)\to S^{\prime
0}_\alpha(K')$ (see the paragraph after Definition~\ref{d1}), constitute an inductive system which will be denoted by $\mathcal S_\alpha$. Let
$\mathcal S^{\mathrm{pr}}_\alpha$ be the restriction of $\mathcal S_\alpha$ to $\mathcal P(\R^k)$, i.e., $\mathcal S^{\mathrm{pr}}_\alpha = \mathcal
S_\alpha^{\mathcal P(\R^k)}$. Let $\theta\colon\mathcal P(\R^k)\to \mathcal K(\R^k)$ be the inclusion mapping. It follows from Lemma~\ref{l1} that
$\theta(\mathcal S^{\mathrm{pr}}_\alpha)$ is canonically isomorphic to $\mathcal S_\alpha$ for $\alpha>1$, while Definition~\ref{dxx3} implies that
$\theta(\mathcal S^{\mathrm{pr}}_1)$ coincides with the inductive system $\mathcal U$ indexed by $\mathcal K(\R^k)$ constituted by the spaces
$\mathcal U(K)$ and the linking mappings $\rho^{\mathcal U}_{K',\,K}$.

Let $\alpha>1$. It easily follows from Theorems~\ref{t2} and~\ref{t3} that $\mathcal S_\alpha$ and $\mathcal U$ are prelocalizable inductive systems.
Since $\mathcal S_\alpha$ is prelocalizable, its restriction $\mathcal S^{\mathrm{pr}}_\alpha$ to $\mathcal P(\R^k)$ is also prelocalizable. The same
is true for $\mathcal S^{\mathrm{pr}}_1$, which is naturally isomorphic to the restriction of $\mathcal U$ to $\mathcal P(\R^k)$.

Let $\mathcal X$ and $\mathcal Y$ be inductive systems indexed by the same set $\Gamma$. A map $l$ from $\mathcal X$ to $\mathcal Y$ is, by
definition, a family of linear maps $l_\gamma\colon \mathcal X(\gamma)\to \mathcal Y(\gamma)$ such that $l_{\gamma'}\rho^{\mathcal
X}_{\gamma\gamma'}=\rho^{\mathcal Y}_{\gamma\gamma'}l_\gamma$ for any $\gamma\leq\gamma'$. If $\lambda\colon\Gamma\to\Delta$ is a nondecreasing
mapping, then every $l\colon\mathcal X\to\mathcal Y$ uniquely determines a map $\lambda(l)\colon \lambda(\mathcal X)\to \lambda(\mathcal Y)$ such
that $\lambda(l)_\delta \rho^{\mathcal X^{\Gamma_\delta}}_\gamma=\rho^{\mathcal Y^{\Gamma_\delta}}_\gamma l_\gamma$ for any $\delta\in \Delta$ and
$\gamma\in \Gamma_\delta$. To reformulate Theorem~\ref{t4} in terms of abstract inductive systems, we introduce the next definition.

\begin{defn}\label{dxx7}
Let $\mathcal X$ and $\mathcal Y$ be inductive systems over a partially ordered set $\Gamma$. A mapping $l\colon \mathcal X\to\mathcal Y$ is called
regular if the following conditions hold:
\begin{enumerate}
\item[(a)] $l_\gamma$ are injective for all $\gamma\in\Gamma$,

\item[(b)] if $\gamma\leq \gamma'$ and $\rho^{\mathcal Y}_{\gamma\gamma'}y=l_{\gamma'}x'$ for some $y\in \mathcal Y(\gamma)$ and $x'\in \mathcal
X(\gamma')$, then there exists $x\in \mathcal X(\gamma)$ such that $y=l_{\gamma}x$ (which implies, in view of the injectivity of $l_{\gamma'}$, that
$x'= \rho^{\mathcal X}_{\gamma\gamma'}x$).
\end{enumerate}
\end{defn}

For any $\alpha\geq 1$, we have $S^0_\alpha(\{0\})=S^0_1(\{0\})$, and applying Lemma~\ref{l2} to $K=\{0\}$, we conclude that $S^0_1(K')$ is dense in
$S^0_\alpha(K')$ for any closed proper cone $K'$. Hence, the mappings $j^\alpha_K\colon S^{\prime 0}_\alpha(K)\to S^{\prime 0}_1(K)$ defined in the
end of Sec.~\ref{s1} are injective for proper $K$. By Theorem~\ref{t4}, if $K\subset K'$ and $\rho^1_{K,K'}\, u_1=j^\alpha_{K'}\, u_2$ for some
$u_1\in S^{\prime 0}_1(K)$ and $u_2\in S^{\prime 0}_\alpha(K')$, then there is a $u\in S^{\prime 0}_\alpha(K)$ such that $u_1 = j^\alpha_{K}\,u$ and
$u_2=\rho^\alpha_{K,K'}\,u$. Hence the mapping $j^\alpha\colon \mathcal S^{\mathrm{pr}}_\alpha\to \mathcal S^{\mathrm{pr}}_1$ determined by
$j^\alpha_K$ is regular. We shall derive Theorem~\ref{t1} from the next algebraic statement.

\begin{thm}\label{t6}
Let $\Gamma$ be a distributive quasi-lattice, $\Delta$ be a partially ordered set, $\lambda\colon\Gamma\to\Delta$ be a nondecreasing mapping, and
$\mathcal X$ and $\mathcal Y$ be prelocalizable inductive systems indexed by $\Gamma$. Then $\lambda(l)$ is regular for any regular $l\colon\mathcal
X\to\mathcal Y$.
\end{thm}

\begin{proof}[Proof of Theorem~$\ref{t1}$]
As above, let $\theta\colon\mathcal P(\R^k)\to \mathcal K(\R^k)$ be the inclusion mapping. Clearly, $\Gamma=\mathcal P(\R^k)$, $\Delta=\mathcal
K(\R^k)$, $\lambda=\theta$, $\mathcal X=\mathcal S^{\mathrm{pr}}_\alpha$, and $\mathcal Y=\mathcal S^{\mathrm{pr}}_1$ satisfy the conditions of
Theorem~\ref{t6} and, therefore, the mapping $\theta(j^\alpha)$ is regular. For $K\in\mathcal K(\R^k)$, let $s_K^\alpha\colon S^{\prime
0}_\alpha(K)\to \varinjlim_{M\in\mathcal P(K)}S^{\prime 0}_\alpha(M)$ be the canonical isomorphism provided by Lemma~\ref{l1} and $s^\alpha\colon
\mathcal S_\alpha\to \theta(\mathcal S^{\mathrm{pr}}_\alpha)$ be the mapping determined by $s^\alpha_K$. Then $e^\alpha=s^\alpha\theta(j^\alpha)$,
where $e^\alpha\colon \mathcal S_\alpha\to \mathcal U$ is induced by the mappings $e^\alpha_K$ defined in the end of Sec.~\ref{s1}. Hence $e^\alpha$
is regular, which implies, in particular, that $\varepsilon^\alpha=e^\alpha_{\R^k}$ is injective. Let $K$ be a carrier cone of $\varepsilon^\alpha
u$. Then we have $e^\alpha_{\R^k} u=\rho^{\mathcal U}_{K,\R^k}\,\tilde u$ for some $\tilde u\in \mathcal U(K)$. Since $e^\alpha$ is regular, there is
$\hat u\in S^{\prime 0}_\alpha(K)$ such that $u=\rho^{\alpha}_{K,\R^k}\,\hat u$ and, therefore, $u$ is carried by $K$. Theorem~\ref{t1} is proved.
\end{proof}

To prove Theorem~\ref{t6}, we shall need to introduce some additional notation. Given an inductive system $\mathcal X$ indexed by a partially ordered
set $\Gamma$ and a subset $I$ of $\Gamma$, we denote by $T^{\mathcal X}_I$ the set of triples $(x,\gamma,\gamma')$ such that $\gamma,\gamma'\in I$,
$\gamma\leq\gamma'$, and $x\in \mathcal X(\gamma)$. If $(x,\gamma,\gamma')\in T^{\mathcal X}_\Gamma$, then we set $\sigma^{\mathcal
X}(x,\gamma,\gamma')= \iota^{\mathcal X}_\gamma x-\iota^{\mathcal X}_{\gamma'}\rho^{\mathcal X}_{\gamma\gamma'}x$ (recall that $\iota^{\mathcal
X}_\gamma$ is the canonical embedding of $\mathcal X(\gamma)$ into $\oplus_{\gamma'\in \Gamma} \mathcal X(\gamma')$). We denote by $N^{\mathcal X}_I$
the subspace of $\oplus_{\gamma'\in \Gamma} \mathcal X(\gamma')$ spanned by all $\sigma^{\mathcal X}(x,\gamma,\gamma')$ with $(x,\gamma,\gamma')\in
T^{\mathcal X}_I$. For $I\subset \Gamma$, we denote by $M^{\mathcal X}_I$ the subspace $\oplus_{\gamma\in I} \mathcal X(\gamma)$ of the space
$\oplus_{\gamma\in \Gamma} \mathcal X(\gamma)$. Obviously, the space $\varinjlim \mathcal X^I$ is isomorphic to $M^{\mathcal X}_I/N^{\mathcal X}_I$.
Let $j^{\mathcal X}_I$ be the canonical surjection from $M^{\mathcal X}_I$ onto $\varinjlim \mathcal X^I$. If $I\subset J\subset \Gamma$, then
\begin{equation}\label{3:1}
\tau^{\mathcal X}_{I,\,J}j^{\mathcal X}_I x=
j^{\mathcal X}_J x,\quad x\in M^{\mathcal X}_I.
\end{equation}

A subset $J$ of a quasi-lattice $\Gamma$ will be called $\wedge$-closed if $\gamma_1\wedge\gamma_2\in J$ for any $\gamma_1,\gamma_2\in J$. If
$J\subset \Gamma$ is a $\wedge$-closed set and $J'$ is finite subset of $J$, then one can find a finite $\wedge$-closed set $J''\subset J$ containing
$J'$ (for instance, the set consisting of infima of all subsets of $J'$ can be taken as $J''$). We say that a subset $I$ of a partially ordered set
$\Gamma$ is hereditary if the relations $\gamma\in I$ and $\gamma'\leq \gamma$ imply that $\gamma'\in I$. Clearly, every hereditary subset of a
quasi-lattice is $\wedge$-closed. The proof of Theorem~\ref{t6} is based on the next lemma.

\begin{lem}\label{l3:1}
Let $\Gamma$ be a distributive quasi-lattice, $\mathcal X$ and $\mathcal Y$ be prelocalizable inductive systems indexed by $\Gamma$, and $l\colon
\mathcal X\to \mathcal Y$ be a regular mapping. Let $I$ be a hereditary subset of $\Gamma$, $J$ be a $\wedge$-closed subset of $\Gamma$, and $L\colon
M^{\mathcal X}_\Gamma\to M^{\mathcal Y}_\Gamma$ be the mapping induced by $l_\gamma$. Then we have
\begin{equation}\label{3:2}
N^{\mathcal Y}_J\cap(M^{\mathcal Y}_I+L(M^{\mathcal X}_\Gamma))\subset N^{\mathcal Y}_I+L(N^{\mathcal X}_J).
\end{equation}
\end{lem}
\begin{proof}
For any $y\in N^{\mathcal Y}_J$, there is a finite $\wedge$-closed set $J'\subset J$ such that $y\in N^{\mathcal Y}_{J'}$. Hence it suffices to
prove~(\ref{3:2}) for finite $J$. For $\gamma\in J$, let $k(\gamma)$ be the cardinality $|J_\gamma|$ of the set $J_\gamma = \{\gamma'\in J\,|\,
\gamma'\geq\gamma\}$. Obviously, $\gamma=\inf J_\gamma$. Therefore, if $\gamma,\gamma'\in J$, $\gamma\ne\gamma'$, and $k(\gamma')\geq k(\gamma)$,
then $J_{\gamma}\ne J_{\gamma'}$ and, hence, $k(\gamma\wedge\gamma')=|J_{\gamma\wedge\gamma'}|\geq |J_\gamma\cup J_{\gamma'}|>|J_\gamma|=k(\gamma)$.
For $n\in \mathbb N$, set $C_n =\{\gamma\in J\,|\, k(\gamma)\geq n\}$. We have $J=C_1\supset C_2\supset\ldots \supset C_{|J|}= \{\tilde \gamma\}$,
where $\tilde \gamma=\inf J$, and $C_n=\varnothing$ for $n>|J|$. We say that $y\in N^{\mathcal Y}_J$ admits a decomposition of order $n$ if there are
a family of vectors $y_{\gamma\gamma'}\in \mathcal Y(\gamma)$ indexed by the set $\{(\gamma,\gamma') : \gamma\in C_n,\, \gamma'\in J\setminus
I,\,\gamma<\gamma'\,\}$ and an element $\tilde y\in N^{\mathcal Y}_I+L(N^{\mathcal X}_J)$ such that\footnote{Here and below, we assume that the sum
of a family of vectors indexed by the empty set is equal to zero.}
\begin{equation}
\label{3:3}
y= \tilde y + \sum_{\gamma\in C_n,\, \gamma'\in
J\setminus I,\,\gamma<\gamma'}\sigma^{\mathcal Y}(y_{\gamma\gamma'},
\gamma,\gamma').
\end{equation}
If $y$ has a decomposition of order $>|J|$, then $y\in N^{\mathcal Y}_I+L(N^{\mathcal X}_J)$. Therefore, the lemma will be proved as soon as we show
that every $y\in N^{\mathcal Y}_J\cap(M^{\mathcal Y}_I+L(M^{\mathcal X}_\Gamma))$ admits a decomposition of order $n$ for any $n\in \mathbb N$. Since
$I$ is hereditary, every $y\in N^{\mathcal Y}_J$ has a decomposition of order $1$, and we have to show that $y$ has a decomposition of order $n+1$
supposing it has a decomposition of form~(\ref{3:3}) of order $n$. For this, it suffices to establish that $\sigma^{\mathcal Y}(y_{\gamma\gamma'},
\gamma,\gamma')$ has a decomposition of order $n+1$ for every $\gamma\in C_n,\, \gamma'\in J\setminus I$ such that $\gamma<\gamma'$ and
$k(\gamma)=n$. Let $\Lambda=\{\beta\in C_n\,|\,\beta<\gamma',\,\beta\ne\gamma\}$. Since $\gamma'\notin I$, the $\gamma'$-component of $\tilde y-y$ is
equal to $l_{\gamma'} x'$ for some $x'\in \mathcal X(\gamma')$ and by (\ref{3:3}), we have
\begin{equation}\label{3:4}
\rho^{\mathcal Y}_{\gamma\gamma'}\,y_{\gamma\gamma'}+
\sum_{\beta\in\Lambda}
\rho^{\mathcal Y}_{\beta\gamma'}\,y_{\beta\gamma'}=l_{\gamma'} x'.
\end{equation}
If $\Lambda=\varnothing$, then the regularity of $l$ implies that $y_{\gamma\gamma'}=l_\gamma x$ and $x'=\rho^{\mathcal X}_{\gamma\gamma'}x$ for some
$x\in\mathcal X(\gamma)$ and, hence, $\sigma^{\mathcal Y}(y_{\gamma\gamma'}, \gamma,\gamma')=L\left(\sigma^{\mathcal X}(x,\gamma,\gamma')\right)$.
Therefore $\sigma^{\mathcal Y}(y_{\gamma\gamma'}, \gamma,\gamma')$ admits decompositions of all orders. Let $\Lambda\ne \varnothing$ and $\tilde
\beta=\sup\Lambda$ ($\tilde \beta$ is well-defined because $\Lambda$ is a finite set bounded above by $\gamma'$; note that $\tilde \beta$ does not
necessarily belong to $J$). Set $z=\sum_{\beta\in \Lambda} \rho^{\mathcal Y}_{\beta\tilde\beta}\,y_{\beta\gamma'}$. By (\ref{3:4}), we have
\begin{equation}\label{3:5}
\rho^{\mathcal Y}_{\gamma\gamma'}\,y_{\gamma\gamma'}+\rho^{\mathcal Y}_{\tilde\beta\gamma'}\,z=l_{\gamma'} x'.
\end{equation}
Rewriting (\ref{3:5}) in the form $\rho^{\mathcal Y}_{\gamma\vee\tilde \beta,\,\gamma'}(\rho^{\mathcal Y}_{\gamma,\,\gamma\vee\tilde \beta}
\,y_{\gamma\gamma'}+\rho^{\mathcal Y}_{\tilde\beta,\,\gamma\vee\tilde \beta}z)=l_{\gamma'} x'$ and using the regularity of $l$, we conclude that
$x'=\rho^{\mathcal X}_{\gamma\vee\tilde \beta,\,\gamma'}x$ for some $x\in \mathcal X(\gamma\vee\tilde\beta)$. By condition~(II) of
Definition~\ref{dxx6}, there exist $\eta\in \mathcal X(\gamma)$ and $\zeta\in \mathcal X(\tilde \beta)$ such that $x=\rho^{\mathcal
X}_{\gamma,\,\gamma\vee\tilde\beta}\,\eta + \rho^{\mathcal X}_{\tilde \beta,\,\gamma\vee\tilde\beta}\,\zeta$. Hence $x'= \rho^{\mathcal
X}_{\gamma\gamma'}\,\eta + \rho^{\mathcal X}_{\tilde \beta\gamma'}\,\zeta$ and in view of~(\ref{3:5}), we obtain
\[
\rho^{\mathcal Y}_{\gamma\gamma'}(y_{\gamma\gamma'}-l_\gamma \eta)+\rho^{\mathcal Y}_{\tilde\beta\gamma'}(z-l_{\tilde \beta}\zeta)=0.
\]
By~(III), there is a $w\in \mathcal Y(\tilde\beta\wedge\gamma)$ such that $y_{\gamma\gamma'}=l_\gamma \eta+\rho^{\mathcal
Y}_{\tilde\beta\wedge\gamma,\,\gamma}\,w$. Because the quasi-lattice $\Gamma$ is distributive, we have $\tilde\beta\wedge\gamma=\sup_{\beta\in
\Lambda} \beta\wedge\gamma$ and by~(II), there is a family $\{w_\beta\}_{\beta\in\Lambda}$ such that $w_\beta\in \mathcal Y(\beta\wedge\gamma)$ and
$w=\sum_{\beta\in \Lambda} \rho^{\mathcal Y}_{\beta\wedge\gamma,\,\tilde\beta\wedge\gamma}\,w_{\beta}$. We thus have $y_{\gamma\gamma'}=l_\gamma
\eta+\sum_{\beta\in \Lambda}\rho^{\mathcal Y}_{\beta\wedge\gamma,\,\gamma}\,w_{\beta}$ and, consequently,
\begin{multline}\label{3:6}
\sigma^{\mathcal Y}(y_{\gamma\gamma'},\gamma,\gamma')= L\left(\sigma^{\mathcal X}(\eta,\gamma,\gamma')\right) + \\ +\sum_{\beta\in \Lambda}
[\sigma^{\mathcal Y}(w_{\beta},\gamma\wedge\beta,\gamma')- \sigma^{\mathcal Y}(w_{\beta},\gamma\wedge\beta,\gamma)].
\end{multline}
If $\gamma\notin I$, then (\ref{3:6}) gives a decomposition of order $n+1$ for $\sigma^{\mathcal Y}(y_{\gamma\gamma'},\gamma,\gamma')$ because
$k(\gamma\wedge\beta)>k(\gamma)=n$ for any $\beta\in \Lambda$. If $\gamma\in I$, then the desired decomposition is obtained by rewriting~(\ref{3:6})
in the form $\sigma^{\mathcal Y}(y_{\gamma\gamma'},\gamma,\gamma')=\tilde w + \sum_{\beta\in \Lambda} \sigma^{\mathcal
Y}(w_{\beta},\gamma\wedge\beta,\gamma')$, where $\tilde w = L\left(\sigma^{\mathcal X}(\eta,\gamma,\gamma')\right) - \sum_{\beta\in \Lambda}
\sigma^{\mathcal Y}(w_{\beta},\gamma\wedge\beta,\gamma)$ belongs to $N^{\mathcal Y}_I+L(N^{\mathcal X}_J)$. The lemma is proved.
\end{proof}

The above proof is similar to that of Lemma~A.1 in~\cite{AIM}, which states that
\begin{equation}\label{lemmaa1}
N^{\mathcal Y}_\Gamma\cap M^{\mathcal X}_I = N^{\mathcal X}_I
\end{equation}
for any prelocalizable inductive system $\mathcal Y$ indexed by a quasi-lattice $\Gamma$ and any hereditary set $I\subset \Gamma$. In fact, it is
easy to derive~(\ref{lemmaa1}) from Lemma~\ref{l3:1}. Indeed, let $\mathcal X$ be the trivial inductive system defined by the relations $\mathcal
X(\gamma)=0$ and $\rho^{\mathcal X}_{\gamma\gamma'}=0$ and let $l\colon\mathcal X\to\mathcal Y$ be such that $l_\gamma=0$ for all $\gamma\in\Gamma$.
In view of condition~(I) of Definition~\ref{dxx6}, $l$ is a regular mapping and applying Lemma~\ref{l3:1} to $J=\Gamma$, we obtain~(\ref{lemmaa1}).
We are now ready to prove Theorem~\ref{t6}.

\begin{proof}[Proof of Theorem~$\ref{t6}$]
Let $L$ be as in Lemma~\ref{l3:1}. For any $\delta\in\Delta$, we have
\begin{equation}\label{3:7}
\lambda(l)_\delta j^{\mathcal X}_{\Gamma_\delta} x = j^{\mathcal Y}_{\Gamma_\delta} Lx,\quad x\in M^{\mathcal X}_{\Gamma_\delta}.
\end{equation}
Let $\xi\in \lambda(\mathcal X)(\delta)$ be such that $\lambda(l)_\delta \xi=0$. Choose $x\in M^{\mathcal X}_{\Gamma_\delta}$ such that
$\xi=j^{\mathcal X}_{\Gamma_\delta} x$. By~(\ref{3:7}), we have $j^{\mathcal Y}_{\Gamma_\delta} Lx=0$ and, therefore, $Lx\in N^{\mathcal
Y}_{\Gamma_\delta}$. Applying Lemma~\ref{l3:1} to $I=\varnothing$ and $J=\Gamma_\delta$, we obtain $Lx\in L(N^{\mathcal X}_{\Gamma_\delta})$. Since
$L$ is injective, it follows that $x\in N^{\mathcal X}_{\Gamma_\delta}$ and, hence, $\xi=j^{\mathcal X}_{\Gamma_\delta} x=0$. Thus,
$\lambda(l)_\delta$ is injective for any $\delta\in\Delta$. Let $\delta\leq \delta'$ and $\eta\in \lambda(\mathcal Y)(\delta)$ and $\xi'\in
\lambda(\mathcal X)(\delta')$ be such that
\begin{equation}\label{3:8}
\rho^{\lambda(\mathcal Y)}_{\delta\delta'} \eta= \lambda(l)_{\delta'} \xi'.
\end{equation}
Let $y\in M^{\mathcal Y}_{\Gamma_{\delta}}$ and $x'\in M^{\mathcal X}_{\Gamma_{\delta'}}$ be such that $\eta=j^{\mathcal Y}_{\Gamma_{\delta}}y$ and
$\xi'=j^{\mathcal X}_{\Gamma_{\delta'}} x'$. It follows from~(\ref{3:1}), (\ref{3:7}), and~(\ref{3:8}) that $j^{\mathcal Y}_{\Gamma_{\delta'}}
Lx'=\tau^{\mathcal Y}_{\Gamma_{\delta},\Gamma_{\delta'}}j^{\mathcal Y}_{\Gamma_{\delta}}y=j^{\mathcal Y}_{\Gamma_{\delta'}}y$. This implies that
$Lx'-y\in N^{\mathcal Y}_{\Gamma_{\delta'}}$. Applying Lemma~\ref{l3:1} to $I=\Gamma_\delta$ and $J=\Gamma_{\delta'}$, we conclude that $Lx'-y=\tilde
y + L\tilde x$, where $\tilde y\in N^{\mathcal Y}_{\Gamma_{\delta}}$ and $\tilde x\in N^{\mathcal X}_{\Gamma_{\delta'}}$. Let $x=x'-\tilde x$. We
have $Lx=y+\tilde y\in M^{\mathcal Y}_{\Gamma_{\delta}}$ and, therefore, $x\in M^{\mathcal X}_{\Gamma_{\delta}}$. Let $\xi= j^{\mathcal
X}_{\Gamma_\delta} x$. Then $\xi\in\lambda(\mathcal X)(\delta)$ and in view of~(\ref{3:7}), we have $\lambda(l)_\delta\xi =j^{\mathcal
Y}_{\Gamma_\delta}(y+\tilde y)=\eta$. The theorem is proved.
\end{proof}

Note that only a part of conditions of Definition~\ref{dxx6} is used in the proofs of Lemma~~\ref{l3:1} and Theorem~\ref{t6}. In fact, it suffices to
assume that $\mathcal X$ satisfies~(II) and $\mathcal Y$ satisfies~(II) and~(III). At the same time, (I) is essential for deriving
formula~(\ref{lemmaa1}) that lies at the basis of the proof of Theorem~\ref{t3} given in~\cite{AIM}.

\end{document}